# The Bilateral Vandermonde Convolution

Martin Erik Horn, University of Potsdam
Am Neuen Palais 10, D - 14469 Potsdam, Germany
E-Mail: marhorn@rz.uni-potsdam.de

**Abstract**

A new and easy way of deriving Gauss's Generalized Hypergeometric Theorem is presented by using the Bilateral Binomial Theorem.

**First Introduction**

Gauss's Generalized Hypergeometric Theorem states that

$$ {}_2H_2[a,b;c,d;1] = \Gamma\begin{bmatrix} c, d, 1-a, 1-b, c+d-a-b-1 \\ c-a, d-a, c-b, d-b \end{bmatrix} . \qquad (1)$$

Slater [1, p. 181-183] (naming it Generalized Gauss Theorem) gives three ways of deducing this theorem. The first method uses the fact that any generalized hypergeometric series can be expressed as a sum of two ordinary hypergeometric series. The second method uses the evaluation at the poles of an imaginary integral. And the third method uses the bilateral analogue of Dougall's theorem.

**Second Introduction**

We live in a world where uncounted and uncountable combinatorial identities like

$$ g(n,x) = \sum_{k=0}^{n} f(\text{Binomial coefficients}, k, x) \qquad (2) $$

exist. But is it really necessary that $k \in \mathbb{N}$ or $k \in \mathbb{Z}$? A simple way of generalizing lots of these identities is to change the summation bilaterally into

$$ h(n,x) = \sum_{k=-\infty}^{\infty} f(\text{Binomial coefficients}, n, k, x) \qquad (3) $$

with $k \in \mathbb{R}$, resulting in a bilateral relation.

For example, the Binomial Theorem

$$ (1+x)^n = \sum_{k=0}^{n} \binom{n}{k} x^k \qquad \text{with } k \in \mathbb{N} \qquad (4) $$

can be generalized into the Bilateral Binomial Theorem

$$ (1+z)^x = \sum_{k=-\infty}^{\infty} \binom{x}{y} z^k = \Gamma\begin{bmatrix} x+1 \\ y+1, x-y+1 \end{bmatrix} \cdot {}_1H_1[y-x; y+1; -z] \qquad (5) $$

with $x, y \in \mathbb{R}$, $k \in \mathbb{R}$, $z \in \mathbb{C}$ and $|z|=1$ which is shown in [3].

**The Bilateral Vandermonde Convolution**

As indicated by Riordan in [3, p. 8], the Vandermonde Convolution formula (first found by Chu several hundreds earlier than Vandermonde),

$$\binom{n}{m} = \sum_{k=0}^{\max(m,p)} \binom{n-p}{m-k}\binom{p}{k} \tag{6}$$

is perhaps the most widely used combinatorial identity. Therefore the bilateral extension of this theorem, which of course will be Gauss's Generalized Hypergeometric Theorem (1), is of some interest. So let's start deriving it by using the Bilateral Binomial Theorem (5)

$$\sum_{k=-\infty}^{\infty} \binom{n}{k} x^k = (1+x)^n \tag{7}$$

$$= (1+x)^{n-p}(1+x)^p \tag{8}$$

$$= \left[\ldots + \binom{n-p}{k-1}x^{k-1} + \binom{n-p}{k}x^k + \binom{n-p}{k+1}x^{k+1} + \ldots\right]$$
$$\cdot \left[\ldots + \binom{p}{l-1}x^{l-1} + \binom{p}{l}x^l + \binom{p}{l+1}x^{l-1} + \ldots\right] \tag{9}$$

$$= \ldots + \binom{n-p}{k-2}\binom{p}{l}x^{k+l-2} + \binom{n-p}{k-1}\binom{p}{l-1}x^{k+l-2} + \binom{n-p}{k}\binom{p}{l-2}x^{k+l-2} + \ldots$$
$$\ldots + \binom{n-p}{k-1}\binom{p}{l}x^{k+l-1} + \binom{n-p}{k}\binom{p}{l-1}x^{k+l-1} + \binom{n-p}{k+1}\binom{p}{l-2}x^{k+l-1} + \ldots \tag{10}$$
$$\ldots + \binom{n-p}{k}\binom{p}{l}x^{k+l} + \binom{n-p}{k+1}\binom{p}{l-1}x^{k+l} + \binom{n-p}{k+2}\binom{p}{l-2}x^{k+l} + \ldots$$

$$= \ldots + \sum_{\substack{m=-\infty \\ m\in\mathbb{Z}}}^{m=\infty} \binom{n-p}{k+m-2}\binom{p}{l-m}x^{k+l-2} + \sum_{\substack{m=-\infty \\ m\in\mathbb{Z}}}^{m=\infty} \binom{n-p}{k+m-1}\binom{p}{l-m}x^{k+l-1}$$
$$+ \sum_{\substack{m=-\infty \\ m\in\mathbb{Z}}}^{m=\infty} \binom{n-p}{k+m}\binom{p}{l-m}x^{k+l} + \ldots \tag{11}$$

$$= \sum_{\substack{r=-\infty \\ r\in\mathbb{Z}}}^{r=\infty}\left(\sum_{\substack{m=-\infty \\ m\in\mathbb{Z}}}^{m=\infty} \binom{n-p}{k+m-r}\binom{p}{l-m}\right)x^{k+l-r} \tag{12}$$

Now let's change the summation index $r$ (please note, that $k, l \in \mathbb{R}$) using

$$r = k + l - K \quad \text{resp.} \quad k + l - r = K \tag{13}$$

resulting in

$$\sum_{k=-\infty}^{\infty}\binom{n}{k}x^k = \sum_{K=-\infty}^{\infty}\binom{n}{k}x^K = \sum_{\substack{K=-\infty \\ K\in\mathbb{R}}}^{K=\infty}\left(\sum_{\substack{m=-\infty \\ m\in\mathbb{Z}}}^{m=\infty}\binom{n-p}{m+K-l}\binom{p}{l-m}\right)x^K \tag{14}$$

Therefore every binomial coefficient can be expressed as

$$\binom{n}{K} = \sum_{\substack{m=-\infty \\ m\in\mathbb{Z}}}^{m=\infty}\binom{n-p}{m+K-l}\binom{p}{l-m} \tag{15}$$

Now it's time to change the second summation index into

$$m = l - M \qquad \text{resp.} \qquad M = l - m \tag{16}$$

This then results in the Bilateral Vandermonde Convolution Formula, also known as Gauss's Generalized Hypergeometric Theorem:

$$\binom{n}{K} = \sum_{\substack{M=-\infty \\ M \in \mathbb{R}}}^{M=\infty} \binom{n-p}{K-M}\binom{p}{M} \tag{17}$$

$$= \Gamma\begin{bmatrix} n-p+1, p+1 \\ K-M+1, n-p-K+M+1, M+1, p-M+1 \end{bmatrix} \cdot {}_2H_2\begin{bmatrix} K-M, M-p \\ n-p-K+M+1, M+1 \end{bmatrix} \tag{18}$$

Here *M* can be chosen at pleasure.

**Outlook**

Using the Bilateral Multinomial Theorem [4] a further generalization of all this can be found.

**Literature**

[1] Lucy Joan Slater: Generalized Hypergeometric Functions,
Cambridge University Press, Cambridge 1966.

[2] John Riordan: Combinatorial Identities, Wiley Series in Probability and Mathematical Statistics, John Wiley & Sons, New York, London, Sydney 1968

[3] Martin Erik Horn: Bilateral Binomial Theorem,
Siam-Problem No. 001-03,
http://www.siam.org/journals/problems/03-001.htm   (24.07.2003)

[4] Martin Erik Horn: Pascal Pyramids, Pascal Hyper-Pyramids and a Bilateral Multinomial Theorem, math.GM/0311035,
http://www.arxiv.org/abs/math.GM/0311035   (04.11.2003)